\documentclass[a4paper, reqno, 10pt]{amsart}
\usepackage{amssymb,latexsym, mathdots,enumerate}

\newtheorem{theorem}{Theorem}

\newtheorem{conjecture}{Conjecture}

\begin{document}

%

\title{On the Number of Rainbow Spanning Trees in Edge-Colored Complete Graphs}

\author{Hung-Lin Fu}
\address{Department of Applied Mathematics \\
National Chiao Tung University \\
Hsinchu 300, Taiwan} \email[Hung-Lin Fu]{hlfu@math.nctu.edu.tw}

\author{Yuan-Hsun Lo}
\address{School of Mathematical Science \\
Xiamen University \\
Xiamen, 361005, PRC} \email[Yuan-Hsun Lo]{yhlo0830@gmail.com}

\author{K. E. Perry}
\address{Department of Mathematics and Statistics \\
Auburn University \\
Auburn, Alabama 36849, United States} \email[K. E. Perry]{kep0024@auburn.edu}

\author{C. A. Rodger}
\address{Department of Mathematics and Statistics \\
Auburn University \\
Auburn, Alabama 36849, United States} \email[C. A. Rodger]{rodgec1@auburn.edu}

\date{\today}

\begin{abstract}

A spanning tree of a properly edge-colored complete graph, $K_n$, is rainbow provided that each of its edges receives a distinct color. 
In 1996, Brualdi and Hollingsworth conjectured that if $K_{2m}$ is properly $(2m-1)$-edge-colored, then the edges of $K_{2m}$ can be partitioned into $m$ rainbow spanning trees except when $m=2$. 
By means of an explicit, constructive approach, in this paper we construct $\lfloor \sqrt{6m+9}/3 \rfloor$ mutually edge-disjoint rainbow spanning trees for any positive value of $m$.
Not only are the rainbow trees produced, but also some structure of each rainbow spanning tree is determined in the process.
This improves upon best constructive result to date in the literature which produces exactly three rainbow trees.
\end{abstract}



\maketitle


\section{Introduction}
A spanning tree $T$ of a connected graph $G$ is an acyclic connected subgraph of $G$ for which $V(T) = V(G)$. A proper $k$-edge-coloring of a graph $G$ is a mapping from $E(G)$ into a set of colors, $\{1,2,...,k\}$, such that adjacent edges of $G$ receive distinct colors. 
The chromatic index $\chi'(G)$ of a graph $G$ is the minimum number $k$ such that $G$ is $k$-edge-colorable. 
It is well known that $\chi'(K_{2m}) = 2m-1$ and thus, if $K_{2m}$ is properly $(2m-1)$-edge-colored, each color appears at every vertex exactly once. All edge-colorings considered in this paper are proper.

A subgraph in an edge-colored graph is said to be rainbow (sometimes called multicolored or poly-chromatic) if all of its edges receive distinct colors. 
Observe that with any $(2m-1)$-edge-coloring of $K_{2m}$, it is not hard to find a rainbow spanning tree by taking the spanning star, $S_v$, with center $v \in V(K_{2m})$. 
Further, $K_{2m}$ has $m(2m-1)$ edges and it is well known that these edges can be partitioned into $m$ spanning trees. 
This led Brualdi and Hollingsworth \cite{Brualdi} to make the following conjecture in 1996. 

\begin{conjecture}[\cite{Brualdi}]\label{conj:A}
If $K_{2m}$ is $(2m-1)$-edge-colored, then the edges of $K_{2m}$ can be partitioned into $m$ rainbow spanning trees except when $m = 2.$
\end{conjecture}

Based on Brualdi and Hollingsworth's concept, the following related conjectures were proposed in 2002.

\begin{conjecture}[\cite{Constantine}, Constantine]\label{conj:B}
$K_{2m}$ can be edge-colored with $2m-1$ colors in such a way that the edges can be partitioned into $m$ isomorphic rainbow spanning trees except when $m = 2.$
\end{conjecture}

Conjecture~\ref{conj:B} was proved to be true by Akbari, Alipour, Fu, and Lo in 2006 \cite{Akbari}.

\begin{conjecture}[\cite{Constantine}, Constantine]\label{conj:C}
If $K_{2m}$ is $(2m-1)$-edge-colored, then the edges of $K_{2m}$ can be partitioned into $m$ isomorphic rainbow spanning trees except when $m = 2.$ 
\end{conjecture}

\begin{conjecture}[\cite{Kaneko}, Kaneko, Kano, Suzuki]\label{conj:D}
Every properly colored $K_n$ contains $\left \lfloor{\frac{n}{2}}\right \rfloor$ edge-disjoint isomorphic rainbow spanning trees.
\end{conjecture}

Concerning Conjecture~\ref{conj:A}, in \cite{Brualdi}, Brualdi and Hollingsworth proved that every properly $(2m-1)$-edge-colored $K_{2m}$ has two edge-disjoint rainbow spanning trees for $m > 2$, and in 2000, Krussel, Marshall, and Verrall \cite{Krussel} improved this result to three spanning trees. Kaneko, Kano, and Suzuki \cite{Kaneko} then improved the previous result slightly by showing that three edge-disjoint rainbow spanning trees exist in any proper edge-coloring of $K_{2m}$. Recently, Horn \cite{Horn} showed that for $m$ sufficiently large there is an $\epsilon > 0$ such that every properly $(2m-1)$-edge-colored $K_{2m}$ has $\epsilon 2m$ edge-disjoint rainbow spanning trees. And a recently submitted paper by Pokrovskiy and Sudakov \cite{Pokro} shows that every properly $(n-1)$-edge-colored $K_n$ has $\frac{n}{9}$ edge-disjoint spanning rainbow trees.

Balogh, Liu, and Montgomery \cite{Balogh} have also recently submitted a result showing that every properly edge-colored $K_n$ contains at least $\frac{n}{10^{12}}$ edge-disjoint rainbow spanning trees. They then use this result to show that any properly $(2m-1)$-edge-colored $K_{2m}$ contains linearly many edge-disjoint rainbow spanning trees.

Relaxing the restriction that the coloring be proper, Akbari and Alipour \cite{Akbari1} were able to show that two edge-disjoint rainbow spanning trees exist in any edge-coloring of $K_{2m}$ with each color class containing $\leq m$ edges. With the same assumption that the coloring be not necessarily proper and each color appears on at most $m$ edges, Carraher, Hartke, and Horn \cite{Carraher} showed that if $m$ is sufficiently large ($m \geq 500,000$) then $K_{2m}$ contains at least $\left\lfloor \frac{m}{500 \log (2m)}\right \rfloor$ edge-disjoint rainbow spanning trees.

Essentially not much was done on Conjecture~\ref{conj:C} until recently. In 2015 Fu and Lo \cite{Fu} proved that three isomorphic rainbow spanning trees exist in any $(2m-1)$-edge-colored $K_{2m}$, $m \geq 14$ and in 2017, Pokrovskiy and Sudakov \cite{Pokro} proved the existence of $10^{-6}n$ edge-disjoint rainbow spanning $t$-spiders in any properly edge-colored $K_n$, $0.0007n \leq t \leq 0.2n$. Note that a $t$-spider is a tree obtained from a star by subdividing $t$ of its edges once.

In this paper, we focus on Conjecture~\ref{conj:A} by proving that in any $(2m-1)$-edge-coloring of $K_{2m}$, $m \geq 1$, there exist at least $\left\lfloor \frac{\sqrt{6m + 9}}{3}\right \rfloor$ mutually edge-disjoint rainbow spanning trees. 
Asymptotically, this is not as good as the bound in \cite{Carraher} or \cite{Horn}, but our result applies to all values of $m$ and it is better until $m$ is extremely large (over $5.7 \times 10^7$ for the bound in \cite{Carraher}).
Instead of using the probabilistic method to prove the result, as was used in \cite{Carraher} and \cite{Horn}, we derive our bound by means of an explicit, constructive approach. 
So, not only do we actually produce the rainbow trees, but also some structure of each rainbow spanning tree is determined in the process. 
It should be noted that the best constructive result (before ours) is the one in the paper by Krussel, Marshall, and Verrall \cite{Krussel} which produces just three rainbow spanning trees, though the recent result by Pokrovskiy and Sudakov is stronger. 

Here is our main result.

\begin{theorem}\label{thm}
Let $K_{2m}$ be a properly $(2m-1)$-edge-colored graph. 
Then there exist $\Omega_m=\left\lfloor \frac{\sqrt{6m + 9}}{3}\right \rfloor$ mutually edge-disjoint rainbow spanning trees, say $T_1,T_2,\ldots,T_{\Omega_m}$, with the following properties.
\begin{enumerate}[(i)]
\item Each tree has a designated distinct root.
\item The root of $T_1$ has degree $(2m-1) - 2(\Omega_m - 1)$ and has at least $(2m-1) - 4(\Omega_m -1)$ adjacent leaves.
\item For $2 \leq i \leq \Omega_m$, the root of $T_i$ has degree $(2m-1) - i - 2(\Omega_m-i)$ and has at least $(2m-1) - 2i - 4(\Omega_m - i)$ adjacent leaves.
\end{enumerate}
\end{theorem}

It is worth mentioning here that the above conjectures will play important roles in certain applications if they are true. 
Notice that a rainbow spanning tree is orthogonal to the $1$-factorization of $K_{2m}$ (induced by any $(2m-1)$-edge-coloring). 
An application of parallelisms of complete designs to population genetics data can be found in \cite{Banks}. 
Parallelisms are also useful in partitioning consecutive positive integers into sets of equal size with equal power sums \cite{Jacroux}. 
In addition, the discussions of applying colored matchings and design parallelisms to parallel computing appeared in \cite{Harary}.

\section{Proof of Theorem~\ref{thm}}
We will use induction on the number of trees to prove this result. 
We can assume $m \geq 5$ since for $1 \leq m \leq 4$, $\Omega_m = 1$ and the spanning star, $S_r$, in which $r \in V(K_{2m})$ and $r$ is joined to every other vertex, is clearly a rainbow spanning tree of $K_{2m}$. 
When the value of $m$ is clear, it will cause no confusion to simply refer to $\Omega_m$ as $\Omega$. 
It is worth noting that the following inductive proof can be used as a recursive construction to create $\Omega$ rainbow edge-disjoint spanning trees, $T_1, T_2, ..., T_{\Omega}$.\\

For $1 \leq \psi \leq \Omega$, the rainbow edge-disjoint spanning trees, $T_1, T_2, ..., T_{\psi}$, are constructed to satisfy  the proposition $f(\psi)$, defined to be the conjunction of the following three degree and structural characteristics:

\begin{equation}
\label{I1}
\text{\parbox{.85\textwidth}{Each tree has a designated distinct root.}}
\end{equation}

\begin{equation}
\label{I2}
\text{\parbox{.85\textwidth}{The root of $T_1$ has degree $(2m-1) - 2(\psi - 1)$ and has at least $(2m-1) - 4(\psi -1)$ adjacent leaves.}}
\end{equation}

\begin{equation}
\label{I3}
\text{\parbox{.85\textwidth}{For $2 \leq i \leq \psi$, The root of $T_i$ has degree $(2m-1) - i - 2(\psi-i)$ and has at least $(2m-1) - 2i - 4(\psi - i)$ adjacent leaves.}}
\end{equation}

In particular, note here that by ($\ref{I3}$), if $\psi > 1$, then the root of $T_2$ has degree $(2m-1) - 2 - 2(\psi-2) = (2m-1) - 2(\psi -1)$ and at least $(2m-1) - 4 - 4(\psi - 2) = (2m-1) - 4(\psi -1)$ adjacent leaves, sharing these characteristics with $T_1$ (as stated in (\ref{I2})).\\

We begin with some necessary notation. All vertices defined in what follows are in $V(K_{2m})$, the given edge-colored complete graph.

The proof proceeds inductively, producing a list of $j$ edge-disjoint rainbow spanning trees from a list of $j-1$ edge-disjoint rainbow spanning trees; so for $1 \leq i \leq j \leq \Omega$, let $T_{i}^{j}$ be the $i^{th}$ rainbow spanning tree of the $j^{th}$ induction step and let $r_i$ be the designated root of $T_i^j$. Notice that $r_i$ is independent of $j$. 


Suppose $T$ is any spanning tree of the complete graph $K_{2m}$ with root $r$ containing vertices $y, v, w,$ and $v'$, where $ry$ and $rv$ are distinct pendant edges in $T$ (so $y \neq v$ are leaves of $T$) and $y, v \notin \{w, v'\}$ (note that $w$ could equal $v'$). Then define $T' = T[r;y,v;w,v']$ to be the new graph formed from $T$ with edges $ry$ and $rv$ removed and edges $yw$ and $vv'$ added. Formally, $T' = T[r;y,v;w,v'] = T - ry - rv + yw + vv'$. We note here that $T'$ is also a spanning tree of $K_{2m}$ because $y$ and $v$ are leaves in $T$, and thus adding edges $yw$ and $vv'$ does not create a cycle in $T'$.

Our inductive strategy will be to assume we have $k-1$ (where $1 < k \leq \Omega$) edge-disjoint rainbow spanning trees with suitable characteristics satisfying proposition $f(k-1)$ that yield properties ($\ref{I1}$), ($\ref{I2}$), and ($\ref{I3}$) with $\psi = k-1$. From those trees we will construct $k$ edge-disjoint rainbow spanning trees with suitable characteristics that allow properties ($\ref{I1}$), ($\ref{I2}$), and ($\ref{I3}$) to be eventually established when $\psi = k$, thus satisfying $f(k)$.

For this construction, given any $T_i^{j-1}$ with root $r_i$ and distinct pendant edges $r_iy_i^j$ and $r_iv_i^j$, we define $T_i^j$ in the following way:

\begin{equation}
\label{origtree}
T_i^j = T_i^{j-1}[r_i; y_i^j, v_i^j; w_i^j, (v')_i^{j}] = T_i^{j-1} - r_{i}y_{i}^j - r_{i}v_{i}^j + y_{i}^jw_{i}^j + v_{i}^j(v')_i^{j}
\end{equation}

\noindent The choice of the vertices defined in ($\ref{origtree}$) will eventually be made precise, based on the discussion which follows.

When the value of $j$ is clear, it will cause no confusion to refer to the vertices $y_i^j, v_i^j; w_i^j, (v')_i^{j}$ by omitting the superscript and instead writing $T_i^j = T_i^{j-1}[r_i; y_i, v_i; w_i, v'_i]$. We now make the following remarks about the definition of $T_i^j$ above. Recall that for $1 \leq i \leq j \leq \Omega$, $r_i$ is independent of $j$, and thus is the root of both $T_i^{j-1}$ and $T_i^j$. The following is easily seen to be true.

\begin{equation}
\label{rainbowspanning}
\text{\parbox{.85\textwidth}{If $\varphi$ is any proper edge-coloring of $K_{2m}$ and $T_{i}^{j-1}$ is a rainbow spanning tree of $K_{2m}$ with root $r_i$ and distinct pendant edges $r_iy_i$ and $r_iv_i$, then $T_{i}^{j}$ as defined in $(\ref{origtree})$ is also a rainbow spanning tree of $K_{2m}$ if $\varphi(r_{i}y_{i}) = \varphi(v_{i}v'_{i})$ and $\varphi(r_{i}v_{i}) = \varphi(y_{i}w_{i}).$}}
\end{equation}

Next, for $1 \leq i \leq j \leq \Omega$, let $L_i^j = \{x \mid xr_i$ is a pendant edge in $T_i^j \}$ (so $x$ is a leaf adjacent to $r_i$ in $T_i^j$). Define

\begin{equation}
\label{leaves_in_TiK-1}
L_{j} = \displaystyle  \bigcap_{i=1}^{j} L_i^j.
\end{equation}

Notice that if $x \in L_j$, then for $1 \leq i \leq j$, $xr_i$ is a pendant edge in $T_i^j$ .\\

We now begin our inductive proof with induction parameter $k$. Specifically we will prove that for $1 \leq k \leq \Omega$ there exist $k$ edge-disjoint rainbow spanning trees, $T_1^k, T_2^k, ..., T_k^k$ satisfying the inductive parameter $f(k)$ (stated below as $f(k-1)$ in the inductive step).\\

Base Step. 
The case $k = 1$ is seen to be true for all properly edge-colored complete graphs, $K_{2m}$, by letting $r_1$ be any vertex in $V(K_{2m})$ and defining $T_1^1 = S_{r_{1}}$, the spanning star with root $r_1$. It is also clear that $S_{r_{1}}$ satisfies $f(1)$ since $r_1$ has degree $2m-1$ and has $2m - 1$ adjacent leaves, as required in ($\ref{I2}$). Property ($\ref{I3}$) is vacuously true.\\

Induction Step.
Suppose that $\varphi$ is a proper edge-coloring of $K_{2m}$ and that for some $k$ with $1 < k \leq \Omega$, $K_{2m}$ contains $k-1$ edge-disjoint rainbow spanning trees, $T_1^{k-1}, T_2^{k-1}, ... , T_{k-1}^{k-1}$, satisfying $f(k-1)$: 

\begin{enumerate}
\item $r_i$ is the root of tree $T_i^{k-1}$ and $r_i \neq r_c$ for $1 \leq i, c < k$, $i \neq c$,
\item $d_{T_1^{k-1}}(r_1) = (2m-1) - 2(k-2)$ and $r_1$ is adjacent to at least $(2m-1) - 4(k-2)$ leaves in $T_1^{k-1}$, and
\item For $2 \leq i \leq k-1$, $d_{T_i^{k-1}}(r_i) = (2m-1) - i - 2(k-1-i)$ and $r_i$ is adjacent to at least $(2m-1) - 2i - 4(k - 1 - i)$ leaves in $T_i^k$.
\end{enumerate}

It thus remains to construct $k$ edge-disjoint rainbow spanning trees satisfying $f(k)$.\\

We note here that $f(k-1)$ and the definition of $L_{k-1}$ in $(\ref{leaves_in_TiK-1})$ guarantee that a lower bound for $\left\vert{L_{k-1}}\right\vert$ can be obtained by starting with a set containing all $2m$ vertices, then removing the $k-1$ roots of $T_1^{k-1}, T_2^{k-2}, ..., T_{k-1}^{k-1}$, the (at most $4(k-2)$) vertices in $V(T_1^{k-1} \backslash \{r_1\})$ which are not leaves adjacent to $r_1$, and for $2 \leq i < k$, the (at most $2i + 4(k-1-i)$) vertices in $V(T_i^{k-1} \backslash \{r_i\})$ which are not leaves adjacent to $r_i$. Formally,

\begin{equation}
\label{Lcardinality}
\begin{split}
\left\vert{L_{k-1}}\right\vert & \geq 2m - (k-1) - 4(k-2) - \sum\limits_{i = 2}^{k-1} (2i + 4(k-1-i))
\\
& = 2m - (k-1) - 4(k-2) - (3k^2 - 11k + 10)
\\
& = 2m - 3k^2 + 6k - 1.
\end{split}
\end{equation}

Knowing $\left\vert{L_{k-1}}\right\vert$ is useful because later (see $(\ref{InequalityToProve})$) we will show that if $\left\vert{L_{k-1}}\right\vert > 6k-5$, then from $T_1^{k-1},$ $T_2^{k-1}, ..., T_{k-1}^{k-1}$ we can construct $k$ rainbow edge-disjoint spanning trees which satisfy proposition $f(k)$. As the reader might expect, it is from here that the bound on $\Omega$ is obtained: it actually follows that since $k \leq \Omega$, $\left\vert{L_{k-1}}\right\vert > 6k-5$.

First, select any two distinct vertices $r_k, w_k^k \in L_{k-1}$; since it will cause no confusion, we will write $w_k$ for $w_k^k$. Set $r_k$ equal to the root of the $k^{th}$ tree, $T_k^k$. Later, $r_kw_k$ will be an edge removed from $T_k^k$. For now, the two special vertices $r_k$ and $w_k$ play a role in the construction of $T_i^k$ from $T_i^{k-1}$ for $1 \leq i < k$. For convenience we explicitly state and observe the following

\begin{equation}
\label{pick_rk_wk}
\text{\parbox{.7\textwidth}{Since $r_k$ and $w_k$ are distinct vertices in $L_{k-1}$ (defined in $(\ref{leaves_in_TiK-1})$), $r_k$ and $w_k$ are leaves adjacent to $r_i$ for $1 \leq i < k$.}}
\end{equation}

For the sake of clarity, having selected $r_k$ and $w_k$, we now discuss how to construct the trees $T_1^k, T_2^k, ... , T_{k-1}^k$ before returning to our discussion of the construction of $T_k^k$ (though in actuality $T_k^k$ is formed recursively as we are constructing $T_1^k, T_2^k, ... , T_{k-1}^k$).\\

For $1 \leq i < k$, we will find suitable vertices $v_i^k, w_i^k,$ and $v_i^{k'}$, which for convenience we refer to as $v_i, w_i,$ and $v'_i$ respectively, and define $T_i^k$ in the following way:

\begin{equation}
\label{treedef}
\text{\parbox{.6\textwidth}{$T_i^k = T_i^{k-1}[r_i;r_k,v_i;w_i,v'_i]$\\ where $\varphi(r_{i}r_{k}) = \varphi(v_{i}v'_{i})$ and $\varphi(r_{i}v_{i}) = \varphi(r_{k}w_{i})$}}
\end{equation}

It is clear by ($\ref{rainbowspanning}$) that for $1 \leq i < k$, since $T_i^{k-1}$ is a rainbow spanning tree of $K_{2m}$, if $v_i$ is chosen so that $v_ir_i$ is a pendant edge in $T_i^{k-1}$ with $v_i \neq r_k$, then $T_i^k$ is also a rainbow spanning tree of $K_{2m}$ (recall from $(\ref{pick_rk_wk})$ that $r_k \in L_{k-1}$, so by ($\ref{leaves_in_TiK-1}$) $r_kr_i$ is a pendant edge in $T_i^{k-1}$). 

\begin{equation}
\label{subset}
\text{\parbox{.7\textwidth}{Further, since $r_k, w_k \in L_{k-1}$, it is clear from ($\ref{treedef}$) that $(1)$ $r_k, v_i \notin L_k$, and $(2)$ all leaves adjacent to $r_i$ in $T_i^k$ are leaves adjacent to $r_i$ in $T_i^{k-1}$. Therefore $\left\vert{L_{k}}\right\vert < \left\vert{L_{k-1}}\right\vert$.}}
\end{equation}

By the induction hypothesis, the trees $T_1^{k-1}, T_2^{k-1}, ..., T_{k-1}^{k-1}$ satisfy $f(k-1)$. We now show that the trees $T_1^k, T_2^k, ..., T_{k-1}^k$ satisfy properties $(\ref{I1})$, $(\ref{I2})$, and $(\ref{I3})$ of $f(k)$. We will construct a $T_k^k$ below in $(\ref{tkkdef})$ and together, $T_1^k, ..., T_{k-1}^k, T_k^k$ will be a collection of trees satisfying $f(k)$.

First, clearly $(\ref{I1})$ is satisfied.  Further, for $1 \leq i < k$, when $T_i^k$ is formed from $T_i^{k-1}$ (see $(\ref{treedef})$), it can easily be seen that the degree of $r_i$ is decreased by $2$ and the number of leaves adjacent to $r_i$ is decreased by at most $4$.

\begin{itemize}
\item[(i.)] $T_1^k$

By our induction hypothesis, we have that $d_{T_1^{k-1}}(r_1) = (2m-1) - 2(k-2)$ and that $r_1$ is adjacent to at least $(2m-1) - 4(k-2)$ leaves in $T_1^{k-1}$. From $(\ref{treedef})$ we have that $d_{T_1^{k}}(r_1) = d_{T_1^{k-1}}(r_1) - 2 = (2m-1) - 2(k-2) - 2 = (2m-1) - 2(k-1)$ and that $r_1$ is adjacent to at least $(2m-1) - 4(k-2) - 4 = (2m-1) - 4(k-1)$ leaves in $T_1^k$. So $(\ref{I2})$ of $f(k)$ is satisfied.

\item[(ii.)] $T_i^k$, $2 \leq i < k$

By our induction hypothesis, we have that $d_{T_i^{k-1}}(r_i) = (2m-1) - i - 2(k-1-i)$ and that $r_i$ is adjacent to at least $(2m-1) - 2i - 4(k - 1 - i)$ leaves in $T_i^k$. From $(\ref{treedef})$ we have that $d_{T_i^{k}}(r_i) = d_{T_i^{k-1}}(r_i) - 2 = (2m-1) - i - 2(k-1-i) - 2 = (2m-1) - i - 2(k-i)$ and that $r_i$ is adjacent to at least $(2m-1) - 2i - 4(k - 1 - i) - 4 = (2m-1) - 2i - 4(k - i)$ leaves in $T_i^k$. So $(\ref{I3})$ of $f(k)$ is satisfied when $2 \leq i < k$. 
\end{itemize}

Lastly, we can observe that once $v_i$ is selected, vertices $w_i$ and $v'_i$ are determined by the required property from $(\ref{treedef})$ that $\varphi(r_{i}r_{k}) = \varphi(v_{i}v'_{i})$ and $\varphi(r_{i}v_{i}) = \varphi(r_{k}w_{i})$.

It remains to ensure that the trees, $T_1^k, T_2^k, ..., T_{k-1}^k$, are all edge-disjoint. This is also proved using the induction hypothesis that $T_1^{k-1}, T_2^{k-1}, ...,$ $T_{k-1}^{k-1}$ are all edge-disjoint.\\

Now, while forming the rainbow edge-disjoint spanning trees, $T_1^k, T_2^k, ...,$ $T_{k-1}^k$, we simultaneously construct the $k^{th}$ rainbow spanning tree, $T_k^k$, from a sequence of inductively defined graphs, $T_k^k(1), T_k^k(2), ... , T_k^k(k) = T_k^k$ where $T_k^k(1) = S_{r_{k}}$ and at the $i^{th}$ induction step, the formation of $T_k^k(i)$ depends on the choice of $v_i$ used in the construction of $T_i^k$. For $2 \leq i \leq k$ define

\begin{equation}
\label{tkk_seq_def}
\text{\parbox{.85\textwidth}{$T_k^k(i) = T_k^k(i-1) - r_kw_i + w_iw'_i$,\\ where $\varphi(w_1w'_1) = \varphi(r_kw_k)$ and $\varphi(w_iw'_i) = \varphi(r_kw_{i-1})$ for $2 \leq i \leq k$.}}
\end{equation}

Note that for $1 \leq i \leq k-1$, the choice of $v_i$ determines $T_k^k(i)$; the formation of $T_k^k(k)$ is dictated by $T_k^k(k-1)$ since $w_k'$ is determined by requiring that $\varphi(w_kw'_k) = \varphi(r_kw_{k-1})$. It is worth explicitly stating that

\begin{equation}
\label{tkkdef}
\text{\parbox{.85\textwidth}{$T_k^k = T_k^k(k) = S_{r_{k}} - r_kw_1 - ... - r_kw_k + w_1w'_1 + ... + w_kw'_k$,\\ where $\varphi(w_1w'_1) = \varphi(r_kw_k)$ and $\varphi(w_cw'_c) = \varphi(r_kw_{c-1})$ for $2 \leq c \leq k$}}
\end{equation}

Observe that $T_k^k$ is a rainbow graph since each edge removed from $S_{r_k}$ is replaced by a corresponding edge of the same color. Also, one can easily see that $T_k^k$ has $2m-1$ edges; $d_{T_k^k}(r_k) = (2m-1) - k$ since $r_k \notin \{w_1', w_2', ..., w_k'\}$; and $r_k$ has at least $(2m-1) - 2k$ adjacent leaves. Therefore, condition $(\ref{I3})$ of $f(k)$ is satisfied. So it remains to show that $T_k^k$ is acyclic and contains no edges in the trees $T_i^k$ for $1 \leq i \leq k-1$.

For future reference, it is worth gathering two observations just made into one:

\begin{equation}
\label{***}
\text{\parbox{.85\textwidth}{For $1 \leq i <k$, once $v_i$ is chosen, $T_i^k$ and $T_k^k(i)$ are completely determined by the constructions described in $(\ref{treedef})$ and $(\ref{tkk_seq_def})$ respectively.}}
\end{equation}

Due to the fact highlighted above in $(\ref{***})$, our strategy will be to select a suitable $v_i$ and construct $T_i^k$ from $T_i^{k-1}$, while simultaneously constructing $T_k^k(i)$ from $T_k^k(i-1)$. In doing so, we restrict the choices for each $v_i$ in order to achieve the following three properties:

\begin{itemize}
\item[(C1)] The edges in $T_a^k, 1 \leq a < i$ do not appear in $T_i^k$,
\item[(C2)] The edges in $T_k^k$ do not appear in $T_i^k$, $1 \leq i < k$, and
\item[(C3)] $T_k^k$ is acyclic
\end{itemize}

To that end, we let 

\begin{equation}
\label{chooseset}
L_{k-1}^* = L_{k-1} \backslash \{r_k, w_k\}
\end{equation}
and let $v_i$ be any vertex for which the following properties are satisfied (so by $(\ref{***})$, this choice completes the formation of $T_i^k$ and $T_k^k(i)$ for $1 \leq i < k$):

\begin{itemize}
\item[(R1)] $v_i \in L_{k-1}^*$,
\item[(R2)] For $1 \leq c < k$, $c \neq i$, $\varphi(v_{i}r_{c}) \neq \varphi(r_{i}r_{k})$,

\item[(R3)] For $1 \leq a < i$, $\varphi(v_{i}r_{i}) \neq \varphi(r_{a}v_{a})$, 
\item[(R4)] For $i < b < k$, $\varphi(v_{i}r_{i}) \neq \varphi(r_{k}r_{b})$, 
\item[(R5)] $\varphi(v_{i}r_{i}) \neq \varphi(r_{k}w_k)$,
\item[(R6)] For $1 \leq a < i$, $\varphi(v_{i}r_{i}) \neq \varphi(r_{k}w'_a)$,
\item[(R7)] For $2 \leq i < k$, $\varphi(v_{i}r_{i}) \neq \varphi(r_{k}\alpha)$,\\ where $\alpha$ is the vertex such that $\varphi(w_{k}\alpha) = \varphi(r_kw_{i-1})$,
\item[(R8)] For $i = 1$ and for $1 \leq c < k$, $\varphi(v_{1}r_{1}) \neq \varphi(r_{k}\alpha)$,\\ for each vertex $\alpha$ incident with the edge of color $\varphi(r_{k}w_{k})$ in $T_c^{k-1}$,
\item[(R9)] For $2 \leq i < k$, $1 \leq a < i$, and for $i \leq b < k$, $\varphi(v_{i}r_{i}) \neq \varphi(r_{k}\alpha)$,\\ for each vertex $\alpha$ incident with the edge of color $\varphi(r_{k}w_{i-1})$ in $T_a^{k}$ and in $T_b^{k-1}$,
\item[(R10)] For $1 \leq i < k$, $\varphi(v_iw_k) \neq \varphi(r_ir_k)$,
\item[(R11)] For $1 \leq d \leq k-2$, $\varphi(v_{k-1}r_{k-1}) \neq \varphi(w_kr_d)$.
\end{itemize}

From the observation in $(\ref{Lcardinality})$, we know that $\left\vert  L_{k-1}^* \right\vert \geq 2m - 3k^2 + 6k - 3$.

As $i$ increases, the number of vertices eliminated in each item increases or is constant, except for (R4) and (R8). However, we observe here that the number eliminated in items (R3) and (R4) is $(i-1) + (k-i-1) = (k-2)$, a constant, and the number eliminated in (R8) and (R9) is $2i$, which is maximized when $i = k-1$. Therefore, an upper bound for the number of vertices eliminated through items (R2) - (R11) as candidates for $v_i$ is achieved when $i = k -1$. 
When $i = k-1$, the number of vertices eliminated by (R2), (R3), ..., (R11) is $(k-2),(k-2), 0, 1, (k-2), 1, 0, 2(k-1), 1, (k-2)$ respectively, the sum of which is $6k - 7$. Now, since the induction hypothesis includes the condition $k \leq \Omega$, we can observe the following.

First, from $f(\Omega)$ and the definition of $L_{\Omega-1}$, we can follow the same steps as we did in $(\ref{Lcardinality})$ to see that $\left\vert{L_{\Omega-1}}\right\vert \geq 2m - 3\Omega^2 + 6\Omega - 1$ and further, that $\left\vert  L_{\Omega-1}^* \right\vert \geq 2m - 3\Omega^2 + 6\Omega - 3$. Now, since by the induction hypothesis $k \leq \Omega$ and by ($\ref{subset}$) and ($\ref{chooseset}$), $\left\vert{L^*_{i-1}}\right\vert > \left\vert{L^*_{i}}\right\vert$ for $2 \leq i \leq k-1$, we have the following by our choice of $\Omega$:

\begin{equation}
\label{InequalityToProve}
\left\vert{L^*_{k-1}}\right\vert \geq \left\vert{L^*_{\Omega-1}}\right\vert 
\geq 2m - 3\Omega^2 + 6\Omega - 3
> 6\Omega - 7
\geq 6k - 7.
\end{equation}

Therefore, since $\left\vert  L_{k-1}^* \right\vert > 6k - 7$, such a vertex $v_i$ meeting the restrictions in (R1) - (R11) exists. The following cases show that this choice of $v_i$ ensures that (C1), (C2), and (C3) hold.


\subsection{Case 1}
\textbf{(C1) Edges in $T_a^k, 1 \leq a < i$ do not appear in $T_i^k$.}\\
First, by the induction hypothesis we know that the trees $T_1^{k-1}, T_2^{k-1}, ...,$ $T_{k-1}^{k-1}$ are all rainbow edge-disjoint and spanning. Inductively, we also assume for some $i$ with $2 \leq i < k$ the trees $T_1^k, T_2^k, ..., T_{i-1}^k$ are edge-disjoint rainbow spanning trees as well. By $(\ref{treedef})$, regardless of the choice of $v_i$, the only edges in $T_i^k$ ($1 \leq i < k$) that are not in $T_i^{k-1}$ are $v_{i}v'_{i}$ and $r_{k}w_{i}$. Thus, if we can prove that the edges in $(E(T_i^{k-1}) \backslash \{r_iv_i, r_ir_k\}) \cup \{v_iv'_i, r_kw_i\}$ are not in $T_a^k$, $1 \leq a < i$, we will have shown that the trees $T_1^k, T_2^k, ..., T_i^k$ are all edge-disjoint rainbow and spanning; so by induction, $T_1^k, T_2^k, ..., T_{k-1}^k$ are edge-disjoint rainbow spanning trees.

To that end, for the remainder of Case 1 suppose that $2 \leq i < k$, $1 \leq a < i$, and $i < b < k$ and define the following sets of edges.

\begin{enumerate}
\item $E_{old}(T_a^k) = \{xy \mid xy \in E(T_a^{k-1}) \cap E(T_a^k)\}$
\item $E_{new}(T_a^k) = E(T_a^k) \backslash E(T_a^{k-1}) =  \{v_{a}v'_{a}$, $r_{k}w_{a}\}$
\item $E_{old}(T_i^k) = \{xy \mid xy \in E(T_i^{k-1}) \cap E(T_i^k)\}$
\item $E_{new}(T_i^k) = E(T_i^k) \backslash E(T_i^{k-1}) = \{v_{i}v'_{i}$, $r_{k}w_{i}\}$
\end{enumerate}

Observe that by $(\ref{treedef})$, $E_{old}(T_a^k) \cap E_{new}(T_a^k) = \emptyset$ and $E(T_a^k) = E_{old}(T_a^k) \cup E_{new}(T_a^k)$. Similarly, $E_{old}(T_i^k) \cap E_{new}(T_i^k) = \emptyset$ and $E(T_i^k) = E_{old}(T_i^k) \cup E_{new}(T_i^k)$.

Since the trees $T_1^k, T_2^k, ..., T_{k-1}^k$ are formed sequentially, it is clearly necessary to prohibit edges $v_{i}v'_{i}$ and $r_kw_i$ from appearing in $T_a^k$. It is also very useful to prohibit edges $v_iv'_i$ and $r_kw_i$ from appearing in $T_b^{k-1}$.



Consequently, when $v_i$ was selected to satisfy (R1) - (R11) it was done in such a way that ensures $E_{new}(T_i^k) \cap (E(T_a^k) \cup E(T_b^{k-1}) = \emptyset$ and $E_{old}(T_i^k) \cap (E(T_a^k) \cup E(T_b^{k-1}) = \emptyset$. To prove this, six cases are considered:

\begin{itemize}
\item[(P1)] $v_iv'_i$, $r_kw_i \notin E_{old}(T_a^k)$,
\item[(P2)] $v_iv'_i$, $r_kw_i \notin E_{new}(T_a^k)$,
\item[(P3)] $v_iv'_i$, $r_kw_i \notin E(T_b^{k-1})$,
\item[(P4)] $E_{old}(T_i^k) \cap E_{old}(T_a^k) = \emptyset$,
\item[(P5)] $E_{old}(T_i^k) \cap E_{new}(T_a^k) = \emptyset$,
\item[(P6)] $E_{old}(T_i^k) \cap E(T_b^{k-1}) = \emptyset$.
\end{itemize}

It is clear that if properties (P1) - (P6) are satisfied, then $T_i^k$ is edge-disjoint from the trees, $T_a^k$ and $T_b^{k-1}$. We consider edges $v_iv'_i$ and $r_kw_i$ in turn for properties (P1) - (P3), then address properties (P4) - (P6).

\subsubsection{Properties (P1) and (P3) for $v_iv'_i$}
\label{p1_viv'i}
\textbf{}\\
Since $E_{old}(T_a^k) \subset E(T_a^{k-1})$, we can prove $v_iv'_i$ is not an edge in $E_{old}(T_a^k)$ or $T_b^{k-1}$ by showing that $v_iv'_i \notin E(T_c^{k-1})$ for $1 \leq c < k$, $c \neq i$.

Recall from (R1) and $(\ref{chooseset})$ that because $v_i \in L^*_{k-1}$, $v_i$ is a leaf adjacent to the root $r_c$ in $T_c^{k-1}$. Therefore, to show that $v_iv'_i \notin E(T_c^{k-1})$, we need only prove that $v'_i \neq r_c$. The following argument shows that (R2) guarantees this property.

Suppose to the contrary that $v_i' = r_c$. Then $v_iv_i' = v_ir_c$ and by ($\ref{treedef})$, $\varphi(v_{i}r_{c}) = \varphi(v_iv'_i) = \varphi(r_{i}r_{k})$, contradicting (R2). It follows that $v_i' \neq r_c$ so $v_iv'_i \notin E(T_c^{k-1})$, as required.

\subsubsection{Property (P2) for $v_iv'_i$}
\label{p2_viv'i}
\textbf{}\\
Recall that $E_{new}(T_a^k) = \{v_av'_a, r_kw_a\}$. Thus, to prove that $v_iv'_i \notin E_{new}(T_a^k)$ for $1 \leq a < i$, we need only show that $v_iv'_i \neq v_av'_a$ and $v_iv'_i \neq r_kw_a$. We consider each in turn.

\begin{itemize}
\item[(i.)] $v_iv'_i \neq v_av'_a$

By $(\ref{treedef})$, we have that $\varphi(v_{i}v'_{i}) = \varphi(r_{i}r_{k})$ and $\varphi(v_{a}v'_{a}) = \varphi(r_{a}r_{k})$. But, by property (1) of $f(\psi)$ when $\psi = k-1$ we know $r_i \neq r_a$ and so $\varphi(r_{i}r_{k}) \neq \varphi(r_{a}r_{k})$. It follows that $\varphi(v_{i}v'_{i}) \neq \varphi(v_{a}v'_{a})$ and, therefore, $v_iv'_i \neq v_av'_a$.

\item[(ii.)] $v_iv'_i \neq r_kw_a$

Assume that $v_iv'_i = r_kw_a$ and recall from $(\ref{chooseset})$ that because $v_i \in L_{k-1}^*$, $v_i \neq r_k$. Therefore, $v_i = w_a$. By ($\ref{treedef}$), $\varphi(v_{i}v'_{i}) = \varphi(r_{k}r_{i})$, so since we are assuming that $v_{i}v'_{i} = r_{k}w_{a}$, clearly $\varphi(r_{k}r_{i}) = \varphi(r_{k}w_{a})$ and so $w_a = r_i = v_i$. But because $v_i \in L_{k-1}^*$, $v_i \neq r_i$ and this is a contradiction.
\end{itemize}

\noindent Combining the above two arguments, it is clear that  $v_iv'_i \notin E_{new}(T_a^k)$, as required.

\subsubsection{Property (P1) for $r_kw_i$}
\label{p1_rkwi}
\textbf{}\\
Recall from $(\ref{leaves_in_TiK-1})$ that $r_k \in L_{k-1}$, so $r_kr_a$ is a pendant edge in $T_a^{k-1}$ with leaf $r_k$. Therefore, from $(\ref{treedef})$ it is clear that $r_kr_a \notin E(T_a^k)$ since it is removed from $T_a^{k-1}$ in forming $T_a^k$. So $r_k$ is not incident with any edges in $E_{old}(T_a^k)$ and thus, $r_kw_i$ cannot be an edge in $E_{old}(T_a^k)$, as required.

\subsubsection{Property (P2) for $r_kw_i$}
\label{p2_rkwi}
\textbf{}\\
Recall that $E_{new}(T_a^k) = \{v_av'_a, r_kw_a\}$. To show that $r_kw_k \notin E_{new}(T_a^k)$, we prove that $r_kw_i \neq r_kw_a$ and $r_kw_i \neq v_av_a'$ for $1 \leq a < i$. We consider each in turn.

\begin{itemize}
\item[(i.)] $r_kw_i \neq r_kw_a$\\
To show that $r_kw_i \neq r_kw_a$, we need only show that $w_i \neq w_a$.

By $(\ref{treedef})$ we have that $\varphi(r_{k}w_{i}) = \varphi(r_{i}v_{i})$ and $\varphi(r_{k}w_{a}) = \varphi(r_{a}v_{a})$. So if $r_kw_i = r_kw_a$, then $\varphi(v_{i}r_{i}) = \varphi(r_{a}v_{a})$, contradicting (R3). Therefore, $r_kw_i \neq r_kw_a$, as required.

\item[(ii.)] $r_kw_i \neq v_av'_a$\\
Assume that $r_kw_i = v_av'_a$. Recall from $(\ref{chooseset})$ that because $v_a \in L_{k-1}^*$, $v_a \neq r_k$. Therefore, $v_a = w_i$. By ($\ref{treedef}$), $\varphi(v_{a}v'_{a}) = \varphi(r_{a}r_{k})$, so since we are assuming that $r_kw_i = v_av'_a$, then $\varphi(r_{k}w_{i}) = \varphi(r_{k}r_{a})$ and it follows that $r_a = w_i = v_a$. But this is a contradiction because $v_a \in L_{k-1}^*$ so by $(\ref{chooseset})$, $v_a \neq r_a$.
\end{itemize}

\noindent Combining the above two arguments, it is clear that  $r_kw_i \notin E_{new}(T_a^k)$, as required.

\subsubsection{Property (P3) for $r_kw_i$}
\label{p3_rkwi}
\textbf{}\\
Recall that by $(\ref{pick_rk_wk})$, because $r_k$ was chosen to be in $L_{k-1}$, $r_k$ is a leaf adjacent to the root of $T_b^{k-1}$, $i < b < k$. Thus, to show $r_kw_i \notin E(T_b^{k-1})$, we need only prove that $w_i \neq r_b$.

By $(\ref{treedef})$, we have that $\varphi(r_{k}w_{i}) = \varphi(v_{i}r_{i})$. So if $w_i = r_b$, then $r_kw_i = r_kr_b$ and $\varphi(v_{i}r_{i}) = \varphi(r_{k}r_{b})$, contradicting (R4). Therefore, $r_kw_i \notin E(T_b^{k-1})$, as required.

\subsubsection{Properties (P4), (P5), and (P6)}
\label{properties_for_EoldTiK}
\textbf{}\\
We consider each property, (P4), (P5), and (P6), in turn.

\begin{itemize}
\item[(i.)] Property (P4)\\
By our induction hypothesis, the trees, $T_1^{k-1}, T_2^{k-1}, ..., T_{k-1}^{k-1}$ are all edge disjoint. So (P4) follows because $E_{old}(T_i^k) \subset E(T_i^{k-1})$ and $E_{old}(T_a^k) \subset E(T_a^{k-1})$.

\item[(ii.)] Property (P5)\\
Since $a < i$, from (P3) (replacing $i$ with $a$), it follows that $\{v_av'_a, r_kw_a\} \cap E(T_c^{k-1}) = \emptyset$, for $a < c < k$. In particular, since $i > a$, it follows that $E_{new}(T_a^k) \cap E(T_i^{k-1}) = \emptyset$. And lastly, since $E_{old}(T_i^k) \subset E(T_i^{k-1})$, we have that $E_{old}(T_i^k) \cap E_{new}(T_a^k) = \emptyset$.

\item[(iii.)] Property (P6)\\
Again, by our induction hypothesis, the trees, $T_1^{k-1}, T_2^{k-1}, ..., T_{k-1}^{k-1}$ are all edge-disjoint. It follows that $E_{old}(T_i^k) \cap E(T_b^{k-1}) = \emptyset$ because $E_{old}(T_i^k) \subset E(T_i^{k-1})$.
\end{itemize}

Therefore, properties (P4) - (P6) hold for $E_{old}(T_i^k)$.\\

\noindent The above Sections $\ref{p1_viv'i} - \ref{properties_for_EoldTiK}$ ensure that properties (P1) - (P6) hold. As stated above, since these six properties hold, the trees $T_1^k, T_2^k, ..., T_{k-1}^k$ are all edge-disjoint and further, from $(\ref{treedef})$, are also rainbow and spanning.

\subsection{Case 2}
\textbf{(C2) Edges in $T_k^k$ do not appear in $T_i^k$.}\\
Recall from $(\ref{tkk_seq_def})$ that $T_k^k$ is defined by a sequence, $T_k^k(1), T_k^k(2), ..., T_k^k(k)$, and from $(\ref{***})$ that at the $i^{th}$ induction step, $T_k^k(i)$ was determined by the choice of $v_i$. For the remainder of Case 2, suppose that $1 \leq i < k$, $1 \leq a < i$, and $i < b < k$.\\


In order to prevent edges in $T_k^k$ from also appearing in $T_i^k$, we will now show that $T_i^k$ has been constructed in such a way that $T_k^k(i)$ and $T_k^k$ satisfy the following properties: 

\begin{itemize}
\item[(P7)] $E(T_k^k(i)) \cap E(T_a^k) = \emptyset$
\item[(P8)] $E(T_k^k(i)) \cap E(T_b^{k-1}) = \{r_kr_b\}$
\item[(P9)] $E(T_k^k(i)) \cap E_{old}(T_i^k) = \emptyset$
\item[(P10)] $E(T_k^k(i)) \cap E_{new}(T_i^k) = \emptyset$
\item[(P11)] $w_kw_k' \notin E(T_i^k)$

\end{itemize}

We note here that by $(\ref{treedef})$, when $T_b^k$ was constructed from $T_b^{k-1}$, edge $r_kr_b$ was removed, so it does not appear in $T_b^k$. Therefore, it is not necessary to prevent $r_kr_b$ from being an edge in $T_k^k(i)$ nor $T_k^k$. 

Proving the above five properties will be done inductively. We show in the base step that $T_k^k(1)$ satisfies properties (P7) - (P10) with $i = 1$, and then show that for $2 \leq i < k$, $T_k^k(i)$ satisfies the same four properties before finally proving property (P11).

The following preliminary result will be useful in proving properties (P7) - (P11).

\subsubsection{Preliminary Result: $w_i \neq w_k$}
\label{wi_neq_wk}
\textbf{}\\
Recall from $(\ref{pick_rk_wk})$ that $w_k \in L_{k-1}$ was selected with $r_k$ before any of the rainbow spanning trees $T_1^{k-1}, T_2^{k-1}, ..., T_{k-1}^{k-1}$ were revised. It will be useful to show that the vertices $w_i \in T_i^k$, $1 \leq i < k$, cannot equal $w_k$. 

From $(\ref{treedef})$, we have that $\varphi(v_{i}r_i) = \varphi(r_{k}w_i)$. So if $w_i = w_k$, then $\varphi(v_{i}r_{i}) = \varphi(r_{k}w_k)$ contradicting (R5). Therefore, $w_i \neq w_k$.

\subsubsection{Base Step: $i = 1$} 
\label{Tkk1}
\textbf{}\\
Observe that for $2 \leq b < k$, $E(S_{r_k}) \cap E(T_b^{k-1}) = \{r_kr_b\}$ and $E(S_{r_k}) \cap E_{old}(T_1^k) = \emptyset$ since by $(\ref{treedef})$, $r_kr_1$ is removed from $T_1^{k-1}$ when forming $T_1^k$. Further, it is clear from $(\ref{tkk_seq_def})$ that the only edge in $T_k^k(1)$ that is not in $S_{r_k}$ is $w_1w'_1$. 

\begin{itemize}

\item[(i.)] (P7)

Since $i = 1$, there do not exist any such trees $T_a^k$ since $1 \leq a < i$ and so property (P7) is vacuously true. 

\item[(ii.)] (P8) and (P9)

First, recall that $E_{old}(T_1^k) \subset E(T_1^{k-1})$. To establish properties (P8) and (P9), we show that $w_1w_1' \notin E(T_c^{k-1})$ for $1 \leq c < k$.

Suppose to the contrary that $w_1w_1' \in E(T_c^{k-1})$. Recall from $(\ref{tkk_seq_def})$ that $\varphi(w_1w'_1) = \varphi(r_kw_k)$. So if $w_1w'_1 \in E(T_c^{k-1})$, then $w_1$ is a vertex incident to the edge of color $\varphi(r_kw_k)$ in $T_c^{k-1}$. But this is impossible since from $(\ref{treedef})$ we have that $\varphi(v_1r_1) = \varphi(r_kw_1)$ and from (R8) that $\varphi(v_1r_1) \neq \varphi(r_k\alpha)$, where $\alpha$ is a vertex incident to the edge of color $\varphi(r_kw_k)$ in $T_c^{k-1}$. Therefore, $w_1w'_1 \notin E(T_c^{k-1})$ and $T_k^k(1)$ satisfies properties (P8) and (P9).



\item[(iii.)] (P10)

Recall that $E_{new}(T_i^k) = \{v_iv'_i, r_kw_i\}$. To establish (P10) for $T_k^k(1)$, we need only show that $w_1w'_1 \neq v_1v'_1$ and $w_1w'_1 \neq r_kw_1$. We consider each in turn. 

\begin{itemize}
\item[(a.)] $w_1w'_1 \neq v_1v'_1$\\
Recall from $(\ref{treedef})$ that $\varphi(v_{1}v'_{1}) = \varphi(r_{k}r_1)$ and from $(\ref{tkk_seq_def})$ that $\varphi(w_{1}w'_{1}) = \varphi(r_{k}w_k)$. So if $w_1w'_1 = v_1v'_1$, then $\varphi(r_{k}w_k) = \varphi(r_{k}r_1)$ and so $w_k = r_1$. But this is not possible because by $(\ref{pick_rk_wk})$ $w_k \in L_{k-1}$ and so $w_k \neq r_1$. Therefore, $w_1w'_1 \neq v_1v'_1$.

\item[(b.)] $w_1w'_1 \neq r_kw_1$\\
Recall from $(\ref{tkk_seq_def})$ that $\varphi(w_{1}w'_{1}) = \varphi(r_{k}w_k)$. So if $w_1w'_1 = r_kw_1$, then $\varphi(r_{k}w_k) = \varphi(r_{k}w_1)$ and so $w_k = w_1$, contradicting the result in Section $\ref{wi_neq_wk}$. Thus, $w_1w'_1 \neq r_kw_1$.

\end{itemize}

Therefore, property (P10) holds for $T_k^k(1)$ and we have established our base step.

\end{itemize}

\subsubsection{Property (P7) for $2 \leq i < k$}
\label{tkki_p7}
\textbf{}\\
From $(\ref{tkk_seq_def})$, it is clear that the only edge in $T_k^k(i)$ that differs from $T_k^k(i-1)$ is $w_iw'_i$. Therefore, since by induction we have that $T_k^k(i-1)$ satisfies (P7), in order to prove property (P7) is satisfied for $T_k^k(i)$, we need only show that $w_iw'_i$ is not an edge in $T_a^k$, $1 \leq a < i$. 

To that end, suppose to the contrary that $w_iw_i' \in E(T_a^{k})$. Recall from $(\ref{tkk_seq_def})$ that $\varphi(w_iw'_i) = \varphi(r_kw_{i-1})$. So if $w_iw'_i \in E(T_a^{k})$, then $w_i$ is a vertex incident to the edge of color $\varphi(r_kw_{i-1})$ in $T_a^{k}$. But this is impossible since from $(\ref{treedef})$ we have that $\varphi(v_ir_i) = \varphi(r_kw_i)$ and from (R9) that $\varphi(v_ir_i) \neq \varphi(r_k\alpha)$, where $\alpha$ is a vertex incident to the edge of color $\varphi(r_kw_{i-1})$ in $T_a^{k}$. Therefore, $w_iw'_i \notin E(T_a^{k})$ and $T_k^k(i)$ satisfies property (P7).


\subsubsection{Properties (P8) and (P9) for $2 \leq i < k$}
\label{tkki_p8_p9}
\textbf{}\\
Observe again that $E_{old}(T_i^k) \subset E(T_i^{k-1})$. As in Section $\ref{tkki_p7}$, to prove properties (P8) and (P9) for $T_k^k(i)$, we can show that $w_iw'_i \notin E(T_d^{k-1})$, $i \leq d < k$.

For $i \leq d < k$, property (R9), which guarantees $\varphi(v_{i}r_{i}) \neq \varphi(r_{k}\alpha)$, where $\alpha$ is a vertex incident to the edge of color $\varphi(r_{k}w_{i-1})$ in $T_d^{k-1}$, ensures $w_iw'_i \notin E(T_d^{k-1})$, thus ensuring that (P8) and (P9) hold for $T_k^k(i)$. The argument has been omitted here due to its similarity to the argument used above for (P7) in Section $\ref{tkki_p7}$.

\subsubsection{Property (P10) for $2 \leq i < k$}
\label{tkki_p10}
\textbf{}\\
To prove (P10) for $T_k^k(i)$, we need only show that $w_iw'_i \neq v_iv'_i$ and $w_iw'_i \neq r_kw_i$. We consider each in turn. 

\begin{itemize}
\item[(i.)] $w_iw'_i \neq v_iv'_i$\\
Recall from $(\ref{treedef})$ that $\varphi(v_{i}v'_{i}) = \varphi(r_{k}r_i)$ and from $(\ref{tkk_seq_def})$ that $\varphi(w_{i}w'_{i}) = \varphi(r_{k}w_{i-1})$. If $w_iw'_i = v_iv'_i$, then $\varphi(r_{k}w_{i-1}) = \varphi(r_{k}r_i)$ and so $w_{i-1} = r_i$. But $r_kr_i \in E(T_i^{k-1})$ and $r_kw_{i-1} \in E(T_{i-1}^k)$; so if $w_{i-1} = r_i$, this contradicts property (P3) in the $i-1^{th}$ induction step, which in particular (i.e. when $b = i$) ensures that $r_kw_{i-1} \notin E(T_i^{k-1})$. Therefore, $w_iw'_i \neq v_iv'_i$, as required.

\item[(ii.)] $w_iw'_i \neq r_kw_i$\\
Recall from $(\ref{tkk_seq_def})$ that $\varphi(w_{i}w'_{i}) = \varphi(r_{k}w_{i-1})$. If $w_iw'_i = r_kw_i$, then $\varphi(r_{k}w_{i-1}) = \varphi(r_{k}w_i)$ and so $w_{i-1} = w_i$. However, this is impossible by the result in Section $\ref{p2_rkwi}$ which, in particular, proved that $r_kw_i \neq r_kw_a$ for $1 \leq a < i$. Thus, $w_iw'_i \neq r_kw_i$.
\end{itemize}

Therefore, property (P10) holds for $T_k^k(i)$, as required.

\subsubsection{Property (P11) for $w_kw'_k$}
\label{tkki_p11}
\textbf{}\\
The above sections of Case 2 ensure that the rainbow spanning trees $T_1^k, T_2^k, ...,$ $T_{k-1}^k$ and the rainbow spanning graph, $T_k^k(k-1)$ are all edge-disjoint. Thus, it remains to show that $T_1^k, T_2^k, ..., T_{k-1}^k$ and $T_k^k$ are all edge-disjoint. As above, recall from $(\ref{tkk_seq_def})$ that the only edge in $T_k^k$ that differs from $T_k^k(k-1)$ is $w_kw'_k$. Therefore, showing property (P11) holds will prove that $T_1^k, T_2^k, ... T_{k-1}^k$ and $T_k^k$ are edge-disjoint.

First, observe from $(\ref{pick_rk_wk})$ that since $w_k \in L_{k-1}$, $w_k$ is a leaf adjacent to the root $r_i$  in $T_i^{k-1}$ for $1 \leq i < k$. So if $w_kw'_k \in E(T_i^k)$, $w_kw'_k = w_ir_k$, $v_iv'_i$, or $w_kr_i$. We consider each in turn.

\begin{itemize}
\item[(i.)] $w_kw'_k \neq w_ir_k$\\
From $(\ref{pick_rk_wk})$ we know that $w_k \neq r_k$. So if $w_kw'_k = w_ir_k$, then $w_k = w_i$, contradicting the preliminary result in Section $\ref{wi_neq_wk}$. Therefore, $w_kw'_k \neq w_ir_k$, as required.

\item[(ii.)] $w_kw'_k \neq v_iv'_i$\\
Recall from $(\ref{chooseset})$ that since $v_i \in L_{k-1}^*$, $v_i \neq w_k$. So if $w_kw'_k = v_iv'_i$, then $w_k = v'_i$. From $(\ref{treedef})$ we know that $\varphi(v_iv'_i) = \varphi(r_ir_k)$, so if $w_k = v'_i$, then $\varphi(v_iw_k) = \varphi(r_ir_k)$, contradicting (R10). Therefore, $w_kw'_k \neq v_iv'_i$, as required.

\item[(iii.)]$w_kw'_k \neq w_kr_i$\\
Recall from $(\ref{tkk_seq_def})$ that $\varphi(w_kw_k') = \varphi(r_kw_{k-1})$ and suppose that $w_kw'_k = w_kr_i$. First observe that $i \neq k-1$ since $r_kw_{k-1} \in E(T_{k-1}^k)$ and we know from $(\ref{pick_rk_wk})$ and Section $\ref{wi_neq_wk}$ that $w_k \neq r_k$ and  $w_k \neq w_{k-1}$.

Now, for $1 \leq i \leq k-2$, if $w_kw'_k = w_kr_i$ then $r_i = w'_k$. But from $(\ref{treedef})$ and $(\ref{tkk_seq_def})$ if $r_i = w'_k$ then $\varphi(w_kw'_k) = \varphi(r_kw_{k-1}) = \varphi(v_{k-1}r_{k-1}) = \varphi(w_kr_i)$, contradicting (R11). Therefore, $w_kw'_k \neq w_kr_i$, as required.

\end{itemize}

It follows that $w_kw'_k \notin E(T_i^k)$, $1 \leq i < k$.\\

\noindent The above Sections $\ref{wi_neq_wk}$ - $\ref{tkki_p11}$ ensure that the trees $T_1^k, T_2^k, ..., T_{k-1}^k$ and the graph $T_k^k$ are all edge-disjoint. Further, from $(\ref{treedef})$ it is clear that $T_1^k, T_2^k, ..., T_{k-1}^k$ are all rainbow spanning trees and from $(\ref{tkkdef})$ that $T_k^k$ is a spanning rainbow graph (since for every leaf, $w_c$, $1 \leq c \leq k$, which is adjacent to $r_k$ and for which $r_kw_c$ is removed from $T_k^k$, there exists $w_c'$ such that the edge $w_cw'_c$ is added to $T_k^k$ and edge $w_dw'_d$ in $T_k^k$ such that $\varphi(w_dw'_d) = \varphi(r_kw_c)$, where $d \equiv c+1 \mod{k}$.)

\subsection{Case 3}
\textbf{(C3) Preventing cycles from appearing in $T_k^k$.}
\label{tkki_case3}

Properties (C1) and (C2) in the previous sections guarantee that the rainbow spanning trees $T_1^k, T_2^k, ..., T_{k-1}^k$ and the rainbow spanning graph $T_k^k$ are all edge-disjoint. Thus, it remains to prove that $T_k^k$ is acyclic and, therefore, a tree. This is proved inductively, showing that for $1 \leq i \leq k$, $T_k^k(i)$ is acyclic. //


For our base step, we let $T_k^k(0) = S_{r_k}$ and observe that this graph is clearly acyclic. 

Recall from $(\ref{tkk_seq_def})$ that for $1 \leq i \leq k$, $T_k^k(i) = T_k^k(i-1) - r_kw_i + w_iw'_i$. Therefore, since by induction we have that $T_k^k(i-1)$ is acyclic, in order to prove $T_k^k(i)$ is acyclic, we need only show that adding $w_iw'_i$ to $T_k^k(i-1) - r_kw_i$ does not create a cycle. Let $T_k^k(i-1)^* = T_k^k(i-1) - r_kw_i$.

Now, from $(\ref{tkk_seq_def})$ observe that all of the edges in $T_k^k(i-1)$ are of the form $r_kx$, $r_kw_a',$ and $w_aw'_a$, where $1 \leq a < i$ and $x \in V(K_{2m}) \backslash (\displaystyle \{ \bigcup_{a=1}^{i-1} w_a, w'_a\} \cup \{r_k\})$. Thus, $w_i \in \{r_k, x, w_a, w'_a\}$. We now show that $w_i = x$ and, further, that since $w_i = x$, $T_k^k(i)$ is acyclic. //

In showing that $w_i = x$, we first consider the case where $1 \leq i < k$ and then the case where $i = k$ before then showing $T_k^k(i)$ is acyclic for $1 \leq i \leq k$.

\begin{itemize}
\item[(i.)] $w_i = x$, $1 \leq i < k$\\
First observe that $w_i \neq r_k$ since $r_kw_i$ is an edge in $T_i^k$. Also, $w_i \neq w_a$ (this property is established by (R3) and was discussed in Section $\ref{p2_rkwi}$). Lastly, recall from $(\ref{treedef})$ that $\varphi(v_{i}r_{i}) = \varphi(r_{k}w_i)$. So if $w_i = w'_a$ then $\varphi(v_{i}r_{i}) = \varphi(r_{k}w'_{a})$, contradicting (R6). Therefore, $w_i \neq w'_a$ and it follows that $w_i = x$.

\item[(ii.)] $w_k = x$\\
Begin by observing that $w_k \neq r_k$ (since by $(\ref{pick_rk_wk})$ $w_k$ and $r_k$ were chosen to be distinct vertices) and, for $1 \leq i < k$, $w_k \neq w_i$ (this property was established by (R5) and discussed in Section $\ref{wi_neq_wk}$). The following argument shows $w_k \neq w_i'$.

First, observe that $w_k \neq w'_1$ since $\varphi(w_1w'_1) = \varphi(r_kw_k)$, so if $w_k = w'_1$ then $w_1 = r_k$, which we know from $(\ref{treedef})$ cannot be the case.

Now, for $2 \leq i < k$, let $\alpha \in V(K_{2m})$ be the vertex such that $\varphi(w_k\alpha) = \varphi(r_kw_{i-1})$ and recall from $(\ref{tkkdef})$ that $\varphi(w_iw'_i) = \varphi(r_kw_{i-1})$. Suppose that $w_k = w'_i$. Then since $\varphi(w_k\alpha) = \varphi(r_kw_{i-1}) = \varphi(w_iw'_i) = \varphi(w_iw_k)$, $\alpha$ must equal $w_i$. But from $(\ref{treedef})$, we have that $\varphi(v_ir_i) = \varphi(r_kw_i)$, so if $w_i = \alpha$ then $\varphi(v_ir_i) = \varphi(r_k\alpha)$, contradicting (R7) which ensures that $\varphi(v_{i}r_{i}) \neq \varphi(r_{k}\alpha)$, where $\alpha$ is the vertex such that $\varphi(w_{k}\alpha) = \varphi(r_kw_{i-1})$. Therefore, $w_k \neq w'_i$, $2 \leq i < k$.

Combining the above arguments, it is clear that $w_k = x$.

\item[(iii.)] $T_k^k(i)$ is acyclic, $1 \leq i \leq k$\\
Observe that since $w_i = x$, $w_i \in V(K_{2m}) \backslash (\displaystyle \{ \bigcup_{a=1}^{i-1} w_a, w'_a\} \cup \{r_k\})$ and $w_i$ is a leaf adjacent to $r_k$ in $T_k^k(i-1)$. Now, in order for $w_iw'_i$ to create a cycle in $T_k^k(i)$, there would have to exist a path from $w_i$ to $w'_i$ in $T_k^k(i-1)^*$. But, as we just observed, $w_i$ is a leaf in $T_k^k(i-1)$ and since $T_k^k(i-1)^* = T_k^k(i-1) - r_kw_i$, $w_i$ is an isolated vertex in $T_k^k(i-1)^*$ so it follows that no such path exists. Therefore, $T_k^k(i)$ is acyclic, as required.

\end{itemize}

The above three arguments show that $T_k^k$ is acyclic, thus completing the proof of the theorem.

\end{document}